\documentclass{amsart}

\theoremstyle{plain}
\newtheorem{thm}{Theorem}[section]
\newtheorem{lem}[thm]{Lemma}
\newtheorem{pro}[thm]{Proposition}
\newtheorem{cor}[thm]{Corollary}
\theoremstyle{definition}
\newtheorem{exa}[thm]{Example}

\newcommand{\Z}{\mathbb{Z}}
\newcommand{\N}{\mathbb{N}}

\DeclareMathOperator{\GG}{G}

\DeclareMathOperator{\GL}{GL}

\DeclareMathOperator{\tr}{tr}
\DeclareMathOperator{\Lie}{Lie}

\DeclareMathOperator{\End}{End}

\DeclareMathOperator{\Rep}{Rep}
\DeclareMathOperator{\Mat}{Mat}

\newcommand{\J}{\mathcal{J}}
\newcommand{\K}{\mathcal{K}}

\newcommand{\pr}{\Pi^\lambda}

\DeclareMathOperator{\dbslash}{/\!\!/}

\begin{document}
\title[Decomposition of Marsden-Weinstein reductions]{Decomposition of
Marsden-Weinstein reductions for representations of quivers}
\author{William Crawley-Boevey}
\address{Department of Pure Mathematics, University of Leeds, Leeds LS2 9JT, UK}
\email{w.crawley-boevey@leeds.ac.uk}
\subjclass{Primary 16G20; Secondary 53D20}

\begin{abstract}
We decompose the Marsden-Weinstein reductions
for the moment map associated to representations
of a quiver. The decomposition involves symmetric
products of deformations of Kleinian singularities,
as well as other terms.
As a corollary we deduce that the Marsden-Weinstein
reductions are irreducible varieties.
\end{abstract}
\maketitle

\section{Introduction}
Let $K$ be an algebraically closed field of characteristic zero,
and let $Q$ be a quiver with vertex set $I$. If $\alpha\in\N^I$,
the space of representations of $Q$ of dimension vector $\alpha$ is
\[
\Rep(Q,\alpha) = 
\bigoplus_{a\in Q}\Mat(\alpha_{h(a)}\times\alpha_{t(a)},K)
\]
where $h(a)$ and $t(a)$ denote the head and tail vertices 
of an arrow $a$. 
The group 
\[
\GG(\alpha)=(\prod_{i\in I}\GL(\alpha_i,K))/K^*
\]
acts by conjugation on $\Rep(Q,\alpha)$ and on its cotangent
bundle, which may be identified with 
$\Rep(\overline{Q},\alpha)$, where $\overline{Q}$ is the
double of $Q$, obtained from $Q$ by adjoining a reverse 
arrow $a^*:j\to i$ for each arrow $a:i\to j$ in $Q$.
There is a corresponding moment map
\[
\mu_\alpha:\Rep(\overline{Q},\alpha)\to\End(\alpha)_0,
\quad
\mu_\alpha(x)_i = 
\sum_{\substack{a\in Q \\ h(a)=i}} x_a x_{a^*}
-
\sum_{\substack{a\in Q \\ t(a)=i}} x_{a^*} x_a
\]
where 
\[
\End(\alpha)_0 = 
\{
\theta\in\bigoplus_{i\in I}\Mat(\alpha_i,K)
\mid
\sum_{i\in I}\tr(\theta_i)=0
\}
\cong
(\Lie \GG(\alpha))^*,
\]
and the \emph{Marsden-Weinstein reductions} 
(or \emph{symplectic reductions}) are the affine 
quotient varieties
\[
N(\lambda,\alpha) = \mu_\alpha^{-1}(\lambda) \dbslash \GG(\alpha),
\]
where $\lambda$ is an element of $K^I$ with 
$\lambda\cdot\alpha=\sum_{i\in I}\lambda_i\alpha_i$
equal to zero, and it is identified with the element of $\End(\alpha)_0$
whose $i$th component is $\lambda_i \mathrm{I}$.
(Although it is possible to equip $N(\lambda,\alpha)$ with
the structure of a scheme, possibly not reduced, we do not do so 
in this paper.)

We studied this situation in a previous paper \cite{CBmm}, to which
we refer for further information. We showed there that 
$\mu_\alpha^{-1}(\lambda)$ and $N(\lambda,\alpha)$ are nonempty 
if and only if $\alpha\in\N R_\lambda^+$, 
the set of sums (including 0) of elements of the set $R_\lambda^+$
of positive roots $\alpha$ with $\lambda\cdot\alpha=0$ 
(using the root system in $\Z^I$ associated to $Q$, see \cite{Kac}).

The elements of $\mu_\alpha^{-1}(\lambda)$
correspond to modules for a certain algebra $\pr$,
the deformed preprojective algebra of \cite{CBH}, 
and the points of $N(\lambda,\alpha)$
correspond to isomorphism classes of 
semisimple $\pr$-modules of dimension $\alpha$.
In \cite{CBmm} we showed that the possible dimension 
vectors of simple $\pr$-modules are the 
elements of the set
\[
\Sigma_\lambda = 
\{
\alpha\in R_\lambda^+ \mid
\text{$p(\alpha) > \sum_{t=1}^r p(\beta^{(t)})$
whenever $r\ge 2$, $\alpha = \sum_{t=1}^r \beta^{(t)}$ 
and $\beta^{(t)}\in R_\lambda^+$}
\}
\]
where $p(\alpha) = 1 - \alpha\cdot\alpha + 
\sum_{a\in Q} \alpha_{t(a)}\alpha_{h(a)}$.
Moreover, we showed that if $\alpha\in\Sigma_\lambda$ 
then $\mu_\alpha^{-1}(\lambda)$ and $N(\lambda,\alpha)$ are
irreducible varieties of dimension 
$\alpha\cdot\alpha-1+2 p(\alpha)$ and 
$2 p(\alpha)$ respectively. 
For general $\alpha\in\N R_\lambda^+$ it seems that 
$\mu_\alpha^{-1}(\lambda)$ may be rather complicated,
but we show here that $N(\lambda,\alpha)$ is well-behaved.
If $X$ is an affine variety we denote by
$S^m X$ the symmetric product of $m$ copies of $X$.
Our main result is as follows.

\begin{thm}\label{t:mainthm}
Any $\alpha\in\N R_\lambda^+$ has a decomposition
$\alpha = \sigma^{(1)}+\dots+\sigma^{(r)}$
as a sum of elements of $\Sigma_\lambda$, with the property
that any other decomposition of $\alpha$ as a sum of elements
of $\Sigma_\lambda$ is a refinement of this decomposition.
Collecting terms and rewriting this decomposition as
$\alpha = \sum_{t=1}^s m_t \sigma^{(t)}$
where $\sigma^{(1)},\dots,\sigma^{(s)}$ 
are distinct and $m_1,\dots,m_s$ are 
positive integers, we have
\[
N(\lambda,\alpha) \cong \prod_{t=1}^s
S^{m_t} N(\lambda,\sigma^{(t)}).
\]
\end{thm}

The first part of the theorem means that if 
$\alpha = \sum_{j=1}^n \beta^{(j)}$ with  
$\beta^{(j)}\in \Sigma_\lambda$, then 
$\sigma^{(t)}=\sum_{j\in P_t} \beta^{(j)}$
for some partition $\bigcup_{t=1}^r P_t$
of $\{1,\dots,n\}$.

Recall that the roots $\beta$ can be divided into three classes: 
the \emph{real roots} which have $p(\beta)=0$, the
\emph{isotropic imaginary roots} which have $p(\beta)=1$, and the
\emph{non-isotropic imaginary roots} which have $p(\beta)>1$.
We have some observations concerning these classes.

\begin{pro}\label{p:rootobs}
(1) If $\beta$ is a real root in 
$\Sigma_\lambda$, then $N(\lambda,\beta)$ is a point.

(2) If $\beta$ is an isotropic imaginary root in
$\Sigma_\lambda$, then it is indivisible (its components
have no common divisor) and $N(\lambda,\beta)$ is isomorphic 
to a deformation of a Kleinian singularity.

(3) If $\beta$ is a non-isotropic imaginary root 
in $\Sigma_\lambda$ then any positive multiple of $\beta$
is also in $\Sigma_\lambda$. 
\end{pro}

It follows from the proposition (or directly from the proof 
of the theorem) that $m_t=1$ whenever $\sigma^{(t)}$ is
a non-isotropic imaginary root. Thus the theorem actually gives
\[
N(\lambda,\alpha) \cong 
\prod_{\substack{t=1 \\ p(\sigma^{(t)})=1}}^s 
S^{m_t} N(\lambda,\sigma^{(t)})
\times 
\prod_{\substack{t=1 \\ p(\sigma^{(t)})>1}}^s 
N(\lambda,\sigma^{(t)}).
\]

\begin{exa}
If $Q$ is an extended Dynkin quiver with vertex set 
$\{0,1,\dots,n\}$ and $\lambda=0$, then
$\Sigma_0 = \{\delta,\epsilon_0,\dots,\epsilon_n\}$ 
where $\delta$ is the minimal positive
imaginary root and $\epsilon_i$ are the coordinate vectors.
Thus the decomposition of $\alpha\in\N^I$ is
\[
\alpha = \underbrace{\delta+\dots+\delta}_m +
\underbrace{\epsilon_0+\dots+\epsilon_0}_{\alpha_0-m\delta_0}
+\dots+
\underbrace{\epsilon_n+\dots+\epsilon_n}_{\alpha_n-m\delta_n},
\]
where $m$ is the largest integer with $m\delta\le\alpha$.
Thus $N(0,\alpha)\cong S^m N(0,\delta)$,
and $N(0,\delta)$ is the Kleinian 
singularity of type $Q$. See for example 
\cite[Theorem 8.10]{CBH}.
\end{exa}

If $\alpha\in\N R_\lambda^+$, we denote by 
$\vert \alpha \vert_\lambda$
the maximum value of $\sum_{i=1}^n p(\beta^{(i)})$ 
over all decompositions $\alpha=\sum_{i=1}^n \beta^{(i)}$
with the $\beta^{(i)}$ in $R_\lambda^+$.
In fact one may assume that all $\beta^{(i)}$ 
are in $\Sigma_\lambda$,
for amongst all decompositions which realize the maximum,
one that has as many terms as possible clearly has this property.
Now by Theorem~\ref{t:mainthm}, any decomposition of $\alpha$ as 
a sum of elements of $\Sigma_\lambda$ is a refinement of
one special decomposition $\alpha=\sum_{t=1}^r \sigma^{(t)}$.
The defining property of $\Sigma_\lambda$ then shows that
the maximum is only achieved by this special decomposition.
In particular
$\vert \alpha \vert_\lambda = \sum_{t=1}^r p(\sigma^{(t)})$.

Recall that $N(\lambda,\alpha)$ classifies the semisimple
$\pr$-modules of dimension $\alpha$. If $X$ is a semisimple
$\pr$-module, one says that $X$ has representation type
\[
(k_1,\beta^{(1)};\dots;k_n,\beta^{(n)})
\]
if it has composition factors of dimensions $\beta^{(i)}$
occuring with multiplicity $k_i$. Now Theorem~\ref{t:mainthm}
and \cite[Theorems 1.3,1.4]{CBmm} have the following 
immediate consequence.

\begin{cor}\label{c:irred}
If $\alpha\in\N R_\lambda^+$, then $N(\lambda,\alpha)$ is an 
irreducible variety of dimension
$2 \vert \alpha \vert_\lambda$.
The general element of $N(\lambda,\alpha)$ has
representation type
\[
(m_1,\sigma^{(1)}; \dots; m_s,\sigma^{(s)}).
\]
\end{cor}

I would like to thank A. Maffei for some useful 
discussions, and in particular for explaining
Lemma \ref{l:dropsimple} to me.
I would like to thank E. Vasserot for pointing
out the error in an earlier version of this
paper which used schemes instead of varieties.
This research was done in Spring 2000 while
visiting first the program on `Noncommutative Algebra' 
at MSRI, and then the Sonderforschungsbereich on
`Discrete Structures in Mathematics' at Bielefeld 
University. I would like to thank my hosts
at both institutions for their hospitality.

\section{Preliminary results}
Let $Q$ be a quiver with vertex set $I$.
We denote by $(-,-)$ the symmetric bilinear form on $\Z^I$,
\[
(\alpha,\beta) = \sum_{i\in I} 2 \alpha_i \beta_i
- \sum_{a\in\overline{Q}} \alpha_{h(a)} \beta_{t(a)}
\]
and by $q(\alpha) = \frac{1}{2}(\alpha,\alpha)$ the
corresponding quadratic form. Thus $p(\alpha)=1-q(\alpha)$.
We denote by $\epsilon_i\in\N^I$ the coordinate vector
for a vertex $i\in I$. 

If $i$ is a loopfree vertex (so $(\epsilon_i,\epsilon_i)=2$) 
there is a reflection $s_i:\Z^I\to\Z^I$ defined by 
$s_i(\alpha)=\alpha-(\alpha,\epsilon_i)\epsilon_i$,
and a dual reflection $r_i:K^I\to K^I$ with 
$r_i(\lambda)_j = \lambda_j-(\epsilon_i,\epsilon_j)\lambda_i$.
The reflection at vertex $i$ is said to be \emph{admissible} 
for the pair $(\lambda,\alpha)$ if $\lambda_i\neq 0$.
In this case it is shown in \cite{CBH} that 
there are reflection functors relating $\pr$-modules of 
dimension $\alpha$ with $\Pi^{r_i(\lambda)}$-modules
of dimension $s_i(\alpha)$. 
Let $\sim$ be the equivalence relation on $K^I\times\Z^I$
generated by $(\lambda,\alpha)\sim(r_i(\lambda),s_i(\lambda))$
whenever the reflection at $i$ is admissible for $(\lambda,\alpha)$.
We say that $(\nu,\beta)$ is obtained from $(\lambda,\alpha)$
by a sequence of admissible reflections if they are in the same
equivalence class.

\begin{lem}\label{l:reflectionhomeo}
If $(\nu,\beta)$ is obtained from $(\lambda,\alpha)$ by a sequence
of admissible reflections then 
$N(\nu,\beta)\cong N(\lambda,\alpha)$.
\end{lem}

\begin{proof}
This follows from \cite[Lemma 2.2]{CBmm}.
\end{proof}

If $p$ is an oriented cycle in $\overline{Q}$ then
for any $\alpha\in\N^I$ there is a trace function
\[
\tr_p : \Rep(\overline{Q},\alpha)\to K, 
x \mapsto \tr(x_{a_1}\dots x_{a_\ell})
\]
where $p=a_1\dots a_\ell$.
It is invariant under the action of $\GG(\alpha)$.

\begin{lem}\label{l:tracesgen}
If $\lambda\in K^I$ and $\alpha\in\N^I$ 
then the ring of invariants
$K[\mu_\alpha^{-1}(\lambda)]^{\GL(\alpha)}$ is 
generated by the restrictions of the trace functions
$\tr_p$ where $p$ runs through the oriented cycles in 
$\overline{Q}$.
\end{lem}

\begin{proof}
By \cite{LBP} the ring of invariants
$K[\Rep(\overline{Q},\alpha)]^{\GG(\alpha)}$
is generated by the $\tr_p$.
Now $\mu_\alpha^{-1}(\lambda)$ is a closed subvariety
of $\Rep(\overline{Q},\alpha)$, so the restriction
map on functions
\[
K[\Rep(\overline{Q},\alpha)]\to K[\mu_\alpha^{-1}(\lambda)]
\]
is surjective.
Since $\GG(\alpha)$ is reductive and the base field $K$
has characteristic zero, there is a Reynolds operator,
and so it remains surjective on taking invariants.
\end{proof}

The following result was pointed out to the author
by A. Maffei in the context of Nakajima's quiver varieties. 
(The proof is our own.)
If $\lambda_i=0$ we denote by $S_i$ the $\pr$-module 
with dimension vector $\epsilon_i$ in which all arrows
are zero.

\begin{lem}\label{l:dropsimple}
If $i$ is a vertex with $\lambda_i=0$ and $(\alpha,\epsilon_i)>0$,
then any representation of $\pr$ of dimension $\alpha$ has $S_i$ 
as a composition factor, and there is an isomorphism
\[
N(\lambda,\alpha-\epsilon_i) \cong N(\lambda,\alpha).
\]
\end{lem}

\begin{proof}
Since $(\alpha,\epsilon_i)>0$ the vertex $i$ must be loopfree.
Now some composition factor must have dimension 
$\beta$ with $(\beta,\epsilon_i) > 0$. Then
$\beta = \epsilon_i$ by \cite[Lemma 7.2]{CBmm}. 
Since there is no loop at vertex $i$, the 
relevant composition factor is isomorphic to $S_i$.
Now because $\lambda_i=0$, the choice of a decomposition
\[
K^{\alpha_i} \cong K^{\alpha_i-1} \oplus K
\]
induces an embedding
\[
\mu_{\alpha-\epsilon_i}^{-1}(\lambda)
\to
\mu_\alpha^{-1}(\lambda)
\]
and hence a map
$\theta:N(\lambda,\alpha-\epsilon_i)\to N(\lambda,\alpha)$
which by the observation above is a bijection.
We want to prove that is is an isomorphism
of varieties. For this it suffices to prove that it
is a closed embedding. That is, that the map
of commutative algebras
\[
\theta^*:
K[\mu_\alpha^{-1}(\lambda)]^{\GG(\alpha)}
\to
K[\mu_{\alpha-\epsilon_i}^{-1}(\lambda)]^{\GG(\alpha-\epsilon_i)}
\]
is surjective.
Now it is easy to see that this map sends the trace function $\tr_p$
for dimension $\alpha$ to the trace function $\tr_p$ for dimension
$\alpha-\epsilon_i$. Thus the assertion follows 
from Lemma~\ref{l:tracesgen}.
\end{proof}

\section{Symmetric products}
Throughout this section $Q$ is an extended Dynkin quiver,
$\delta$ is its minimal positive imaginary root,
and $\lambda\in K^I$ satisfies $\lambda\cdot\delta=0$.
We choose an extending vertex 0 for $Q$, which
means that $\delta_0=1$. 

We say that an element of the set $\N R_\lambda^+$ is 
\emph{indecomposable} if it is nonzero and it cannot 
be written as a sum of two nonzero elements of this set.

\begin{lem}
The elements of $\Sigma_\lambda$ are $\delta$ and the indecomposable
elements of $\N R_\lambda^+$. All elements are $\le\delta$.
\end{lem}

\begin{proof}
Clearly any real root $\alpha$
in $\Sigma_\lambda$ must be indecomposable since $p(\alpha)=0$.
Conversely, by \cite[Lemma 5.5]{CBmm} any 
indecomposable element is in $\Sigma_\lambda$. 
If $\alpha\in\Sigma_\lambda\setminus\{\delta\}$ is not 
$\le\delta$ then $\alpha-\delta$ is a root with some positive
component, hence a positive root. But 
$\alpha = \delta+(\alpha-\delta)$, contradicting
indecomposability.
\end{proof}

\begin{lem}\label{l:deltasum}
Any decomposition of $m\delta$ as a
sum of elements of $\Sigma_\lambda$ is a refinement
of the decomposition 
\[
m\delta=\underbrace{\delta+\dots +\delta}_m .
\]
\end{lem}

\begin{proof}
Say $\alpha^{(1)},\dots,\alpha^{(q)}$ are elements of $\Sigma_\lambda$
with $\sum_{t=1}^r \alpha^{(t)}=m\delta$.
By induction it suffices to find a subset $P$ of $\{1,\dots,q\}$ with
$\sum_{t\in P}\alpha^{(t)}=\delta$.
We prove this by another induction: if $P$ is
a subset for which the sum is a root $\beta<\delta$, we
show how to enlarge $P$ so that the sum is a root $\le\delta$.
Now $(\delta,\beta)=0$ and $(\beta,\beta)=2$, so
$(\beta,\sum_{t\notin P}\alpha^{(t)})=-2$.
Thus $(\beta,\alpha^{(s)})\le -1$ for some $s\notin P$.
Clearly $\alpha^{(s)}\neq\delta$, so
\[
q(\beta+\alpha^{(s)}) = q(\beta) + q(\alpha^{(s)}) +
(\beta,\alpha^{(s)}) \le 1+1-1 = 1,
\]
so 
$\beta+\alpha^{(s)} = \sum_{t\in P\cup\{s\}}\alpha^{(t)}$ is a root. 
Moreover $\beta+\alpha^{(s)}\le \delta$, for otherwise
$\gamma = \beta+\alpha^{(s)} - \delta$ is a root
(since $q(\gamma)\le 1$) with some positive component,
hence a positive root.
But then $\alpha^{(s)} = \gamma + (\delta-\beta)$, a sum
of elements of $R_\lambda^+$, which contradicts the fact that
$\alpha^{(s)}\in\Sigma_\lambda$.
\end{proof}

\begin{lem}\label{l:througho}
$K[\mu_{\delta}^{-1}(\lambda)]^{\GL(\delta)}$
is generated by the trace functions for paths
in $\overline{Q}$ which start and end at the
extending vertex 0.
\end{lem}

\begin{proof}
Since $\delta_0=1$, the trace function $\tr_p$ for
a path which starts and ends at $0$ involves the trace
of a $1\times 1$ matrix, which is just the
unique entry of the matrix. The assertion thus
follows from \cite[Corollary 8.11]{CBH}. 
\end{proof}

If $X$ is an affine variety, we 
write $S^m X$ for its $m$th \emph{symmetric product}, the
affine variety 
$(X\times\dots\times X)/S_m$. 
Writing $T^m A$ for the $m$th tensor power of an algebra $A$,
we have 
\[
K[S^m X] = (T^m K[X])^{S_m}.
\]

\begin{thm}\label{t:sympower}
The direct sum map 
\[
\prod_{j=1}^m \mu_\delta^{-1}(\lambda) \to \mu_{m\delta}^{-1}(\lambda)
\]
induces an isomorphism
\[
f:S^m N(\lambda,\delta) \to N(\lambda,m\delta) 
\]
\end{thm}

\begin{proof}
By Lemma~\ref{l:deltasum} we know that $f$ is surjective.
Thus it suffices to prove
that it is a closed embedding, that is, that the map on functions
\[
f^*: K[\mu_{k\delta}^{-1}(\lambda)]^{\GL(k\delta)} 
\to \bigl(T^k K[\mu_\delta^{-1}(\lambda)]^{\GL(\delta)}\bigr)^{S_k}
\]
is surjective. 

By Lemma \ref{l:througho} the ring 
$K[\mu_\delta^{-1}(\lambda)]^{\GL(\delta)}$ is 
generated by the trace functions $\tr_p$ for $p$
a path in $\overline{Q}$ starting and ending at 0.
Since the ring is finitely generated, a finite number
of paths $p_1,\dots,p_N$ suffices.

For $1\le j\le m$ let $\pi_j$ be the projection from the product of
$m$ copies of $N(\lambda,\delta)$ onto the $j$th
factor. Thus the coordinate ring of this product is generated by
elements $\tr_{p_i}\circ\pi_j$.

There is a surjective map from the polynomial ring 
$K[x_{ij} : 1\le i\le N, 1\le j\le m]$ to
$T^m (K[\mu_\delta^{-1}(\lambda)]^{\GL(\delta)})$ 
sending $x_{ij}$ to $\tr_{p_i}\circ\pi_j$.
This induces a surjective map 
\[
K[x_{ij}]^{S_m} \to 
\bigl(T^m K[\mu_\delta^{-1}(\lambda)]^{\GL(\delta)}\bigr)^{S_m}
\]
Now by Lemma \ref{l:symminv} below, $K[x_{ij}]^{S_m}$ is 
generated by the power sums
\[
s_{r_1,\dots,r_N} = \sum_j x_{1j}^{r_1} \dots x_{Nj}^{r_N}.
\]
Thus $\bigl(T^m K[\mu_\delta^{-1}(\lambda)]^{\GL(\delta)}\bigr)^{S_m}$
is generated by the elements
\[
s'_{r_1,\dots,r_N} = 
\sum_j (\tr_{p_1}\circ\pi_j)^{r_1} \dots (\tr_{p_N}\circ\pi_j)^{r_N}
= \sum_j (\tr_{p_1}^{r_1} \dots \tr_{p_N}^{r_N})\circ\pi_j.
\]
Since $\delta_0=1$ we have $\tr_p \tr_q = \tr_{pq}$ 
for any paths $p,q$ which
start and end at 0, so $\tr_{p_1}^{r_1} \dots \tr_{p_N}^{r_N} = \tr_p$ 
where $p$ is the path $p_1^{r_1} \dots p_N^{r_N}$. Thus
\[
s'_{r_1,\dots,r_N} = \sum_j \tr_p \circ \pi_j.
\]
This shows that $s'_{r_1,\dots,r_N}$ is the image 
under $f^*$ of the trace function $\tr_p$ for
$\mu_{m\delta}^{-1}(\lambda)$.
Thus the image of $f^*$ contains a set of generators, so $f^*$ is
surjective, as required.
\end{proof}

\begin{lem}\label{l:symminv}
If $S_m$ acts on the polynomial ring 
$K[x_{ij} : 1\le i\le N, 1\le j\le m]$
by permuting the $x_{ij}$ for each $i$, 
then the ring of invariants is generated by
the power sums
\[
s_{r_1,\dots,r_N} = \sum_j x_{1j}^{r_1} \dots x_{Nj}^{r_N}.
\]
($r_1,\dots,r_N\ge 0$).
\end{lem}

\begin{proof}
By \cite[Chapter II, Section 3]{W} the ring of invariants is
generated by polarizations of the elementary symmetric 
polynomials, so by elements of the form
\[
\phi_{i_1,i_2,\dots,i_k} = \sum x_{i_1 j_1} x_{i_2 j_2} \dots x_{i_k j_k}
\]
where the sum is over all distinct 
$j_1,j_2,\dots,j_k$ in the range $1$ to $m$.
Now the elementary symmetric polynomials can be expressed as
polynomials in the power sums by Newton's formulae, 
and on polarizing this expresses 
$\phi_{i_1,i_2,\dots,i_k}$ as a polynomial 
in the $s_{r_1,\dots,r_N}$.
For example polarizing the formula
\[
\sum_{j<k<\ell} z_j z_k z_\ell
=
\frac{1}{6} \left( 
\bigl(\sum_j z_j \bigr)^3 
- 3 \bigl(\sum_j z_j \bigr) \bigl(\sum_j z_j^2 \bigr) 
+ 2 \sum_j z_j^3 
\right)
\]
with respect to the sets of variables $x_{i_1,j}$, 
$x_{i_2,j}$ and $x_{i_3,j}$ gives
\[
\begin{split}
\phi_{i_1,i_2,i_3} = 
& \bigl(\sum_j x_{i_1,j}\bigr) 
\bigl(\sum_j x_{i_2,j}\bigr) 
\bigl(\sum_j x_{i_3,j}\bigr)
- \bigl(\sum_j x_{i_1,j}\bigr)
\bigl(\sum_j x_{i_2,j} x_{i_3,j}\bigr) 
\\
& - \bigl(\sum_j x_{i_2,j}\bigr)
\bigl(\sum_j x_{i_1,j} x_{i_3,j}\bigr)
- \bigl(\sum_j x_{i_3,j}\bigr)
\bigl(\sum_j x_{i_1,j} x_{i_2,j}\bigr)
\\
& + 2 \sum_j x_{i_1,j} x_{i_2,j} x_{i_3,j},
\end{split}
\]
and all sums on the right hand side are of the 
form $s_{r_1,\dots,r_N}$ for suitable $r_1,\dots,r_N$.
\end{proof}

\section{Adding a vertex to an extended Dynkin quiver}\label{s:maincase}
In this section let $Q'$ be an extended
Dynkin quiver, let $k$ be an extending vertex for $Q'$, 
and let $Q$ be a quiver obtained from $Q'$ by adjoining
one vertex $j$ and one arrow joining $j$ to $k$. Let $I$ be the
vertex set of $Q$ and let $\delta\in \N^I$ be the minimal positive
imaginary root for $Q'$.

For any $\alpha\in \Z^I$ we define 
$\alpha' = \alpha - \alpha_j \epsilon_j$.
Thus $\alpha'_j=0$ and $\alpha'_i=\alpha_i$ for $i\neq j$.
One can think of $\alpha'$ as the restriction of $\alpha$
to $Q'$.

Throughout this section
we assume that $\lambda\in K^I$ satisfies
$\lambda\cdot\delta = \lambda_j = 0$.
We prove the following result
which is used in the next section.

\begin{pro}\label{l:maincase}
If $\alpha\in\Sigma_\lambda$, $\alpha_j=1$ and
$m\delta-\alpha'\in\N R_\lambda^+$ for some $m\ge 0$ then
$\alpha=\epsilon_j$.
\end{pro}

An example shows the necessity of
the hypothesis that $m\delta-\alpha'\in\N R_\lambda^+$.

\begin{exa}
Let $Q$ be the quiver
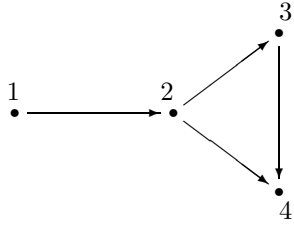
\[
\begin{picture}(120,80)
\put(10,40){\circle*{3}}
\put(70,40){\circle*{3}}
\put(110,10){\circle*{3}}
\put(110,70){\circle*{3}}
\put(15,40){\vector(1,0){50}}
\put(110,65){\vector(0,-1){50}}
\put(74,43){\vector(4,3){32}}
\put(74,37){\vector(4,-3){32}}
\put(7,45){1}
\put(65,45){2}
\put(110,75){3}
\put(110,0){4}
\end{picture}
\]
with vertex set $\{1,2,3,4\}$, so $j=1$, $k=2$,
$Q'$ is of type $\tilde{A}_2$ and
$\delta=(0,1,1,1)$.
If $\lambda=(0,1,-2,1)$ then 
$\alpha=(1,3,2,1) \in\Sigma_\lambda$
since by admissible reflections
at the indicated vertices the pair $(\lambda,\alpha)$
transforms as
\[
\begin{split}
((0,1,-2,1),(1,3,2,1))
\stackrel{2}\sim
((1,-1,-1,2),(1,1,2,1))
\stackrel{3}\sim
((1,-2,1,1),(1,1,0,1))
\\
\stackrel{4}\sim
((1,-1,2,-1),(1,1,0,0))
\stackrel{2}\sim
((0,1,1,-2),(1,0,0,0))
\end{split}
\]
and it is clear that 
$(1,0,0,0)\in \Sigma_{(0,1,1,-2)}$.
However, it is easy to see that there 
is no $m$ with $m\delta-\alpha'=(0,m-3,m-2,m-1)$
in $\N R_\lambda^+$.
\end{exa}

Before proving the proposition we need some lemmas.
Observe that for a vertex $i\notin\{j,k\}$ 
we have $r_i(\lambda)\cdot\delta = 0$ and $r_i(\lambda)_j = 0$.
On the other hand $r_k(\lambda)\cdot\delta = 0$,
but we may have $r_k(\lambda)_j \neq 0$.

\begin{lem}\label{l:gettocoord}
If $\alpha\in\Sigma_\lambda$ and $\alpha_j=1$,
then by a sequence of admissible reflections
at vertices $\neq j$ one can send $(\lambda,\alpha)$ 
to $(\nu,\epsilon_j)$ for some $\nu$.
\end{lem}

\begin{proof}
We consider the pairs $(\nu,\beta)$ which
can be obtained from $(\lambda,\alpha)$
by a sequence of such admissible reflections.
Always $\beta$ is positive, since it is in $\Sigma_{\nu}$
by \cite[Lemma 5.2]{CBmm}. Thus we can choose a pair
$(\nu,\beta)$ with $\beta$ minimal.
Clearly we have $\beta_j=1$.
For a contradiction, suppose that $\beta'\neq 0$.

Since $\delta$ is unchanged by these reflections, we have
$\nu\cdot\delta=\lambda\cdot\delta = 0$.
Also, for each vertex $i\neq j$ we have 
$(\beta,\epsilon_i)\le 0$, for either there is a
loop at $i$, in which case it is automatic, or
$\nu_i=0$, in which case it follows from 
\cite[Lemma 7.2]{CBmm}, or there is an admissible
reflection at $i$, and it follows from the minimality
of $\beta$.
We deduce that $(\beta',\epsilon_i)\le 0$
for $i\notin\{j,k\}$,
and $(\beta',\epsilon_k)\le 1$.

Suppose first that $(\beta',\epsilon_k)=1$. Then
\[
0 = (\beta',\delta) = \sum_{i\neq j} (\beta',\epsilon_i) \delta_i
=
1 + \sum_{i\notin\{j,k\}} (\beta',\epsilon_i) \delta_i,
\]
and all terms in the second sum are $\le 0$. Thus exactly one of
the terms is $-1$, and all others are zero. That is,
there is a vertex $s\neq k$ in $Q'$ with $\delta_s=1$ and
$(\beta',\epsilon_s)=-1$, and $(\beta',\epsilon_i)=0$
for all vertices $i\notin\{k,s\}$ in $Q'$. This is
impossible by \cite[Lemma 8.8]{CBmm}.

Thus $(\beta',\epsilon_k)\le 0$. 
It follows that $(\beta',\beta')\le 0$, so since $Q'$ is
extended Dynkin we have $\beta' = m\delta$ for some $m>0$.
Now the decomposition $\beta = \epsilon_j + \delta+\dots+\delta$
is easily seen to satisfy
\[
p(\beta) = 1 - q(\epsilon_j+m\delta) 
= -m(\epsilon_j,\delta) = m = 
p(\epsilon_j) + p(\delta)+\dots+p(\delta).
\]
We have seen that $\delta\in R_\nu^+$. 
Also
$\nu_j = \nu\cdot\epsilon_j = \nu\cdot\beta = \lambda\cdot\alpha=0$
since $\alpha\in\Sigma_\lambda$, so that $\epsilon_j\in R_\nu^+$.
This contradicts the fact that $\beta\in\Sigma_{\nu}$.

Thus $\beta'=0$, as required.
\end{proof}

\begin{lem}\label{l:rootineq}
If $\alpha\in\Sigma_\lambda$ and $\alpha_j=1$ 
then $\gamma_k-1 \le (\alpha',\gamma) \le \gamma_k$
for any $\gamma\in R_\lambda^+$ with $\gamma < \delta$. 
\end{lem}

\begin{proof}
Some sequence of admissible reflections at vertices $\neq j$
sends $(\lambda,\alpha)$ to $(\nu,\epsilon_j)$.
If $\gamma\in R_\lambda^+$ and $\gamma_j=0$ 
then by \cite[Lemma 5.2]{CBmm} the reflections
send it to a positive root $\beta$, still with $\beta_j=0$.
Thus $(\alpha,\gamma) = (\epsilon_j,\beta) \le 0$,
and so 
\[
(\alpha',\gamma) = (\alpha,\gamma)-(\epsilon_j,\gamma) \le
0-(-\gamma_k) = \gamma_k,
\]
which is one of the inequalities. The other one is obtained
by replacing $\gamma$ with $\delta-\gamma\in R_\lambda^+$.
\end{proof}

Choose a total ordering $\prec$ on $K$ as in
\cite[Section 7]{CBH}.
Let $Q''$ be the Dynkin quiver obtained
from $Q'$ by deleting the vertex $k$.
Let $I''$ be the vertex set of $Q''$.
Recall that a vector $\mu\in K^{I''}$
is said to be \emph{dominant} if $\mu_i \succeq 0$
for all $i\in I''$.

\begin{lem}\label{l:makedominant}
By a sequence of admissible reflections at vertices in $I''$
one can send $(\lambda,\alpha)$ to a pair $(\xi,\beta)$
where $\xi$ is a vector whose restriction to $I''$ is dominant.
\end{lem}

\begin{proof}
Apply \cite[Lemma 7.2]{CBH} to $Q''$, and then consider
the sequence of reflections as reflections for $Q$.
Of course non-admissible
reflections can be omitted, for if $\xi\in K^I$ 
and $\xi_i=0$ then $r_i(\xi)=\xi$.
\end{proof}

\begin{lem}\label{l:dominantcase}
If the restriction of $\lambda$ to $I''$
is dominant, and if $\gamma\in\N R_\lambda^+$ 
has $\gamma_j=0$, then there is some $r\ge 0$
with $\gamma_i=r\delta_i$ for all vertices
$i$ with $\lambda_i\neq 0$.
\end{lem}

\begin{proof}
Any indecomposable element of $\N R_\lambda^+$
which vanishes at $j$ is $\le \delta$, so it 
suffices to prove that if $\gamma\in\N^I$ is
a vector with $\gamma\le\delta$ and 
$\lambda\cdot\gamma=0$, then either $\gamma_i=0$
for all $i$ with $\lambda_i\neq 0$,
or $\gamma_i=\delta_i$ for all $i$ with 
$\lambda_i\neq 0$.

Since $k$ is an extending vertex for $Q'$
we have $\delta_k=1$, and so by replacing
$\gamma$ by $\delta-\gamma$ if necessary,
we may assume that $\gamma_k=0$.

Now the equality $\lambda\cdot\gamma=0$
implies that $\sum_{i\in I''} \gamma_i \lambda_i=0$.
By the dominance condition it follows that $\gamma_i=0$
for any vertex $i\in I''$ with $\lambda_i\neq 0$.
\end{proof}

\begin{proof}[Proof of Proposition \ref{l:maincase}]
First suppose that $\lambda=0$. If $\alpha\neq\epsilon_j$
then the expression for $\alpha$ as a sum of coordinate
vectors is a non-trivial decomposition into elements
of $R_\lambda^+$. Since $p(\alpha)=0$ by Lemma~\ref{l:gettocoord},
this contradicts the fact that $\alpha\in\Sigma_\lambda$.

Thus we may suppose that $\lambda\neq 0$.
Replacing $(\lambda,\alpha)$ by the pair $(\xi,\beta)$ of 
Lemma~\ref{l:makedominant}, we may assume that
the restriction of $\lambda$ to $I''$ is dominant. 
Observe that the reflections involved,
at vertices in $I''$, can change $\alpha$, but they do
not affect the dimension vectors $\epsilon_j$ and $\delta$.
The standing hypotheses on $\lambda$ still hold, as 
do the hypotheses of the proposition
by \cite[Lemma 5.2]{CBmm}.

Now the restriction of $\lambda$ to $I''$ is non-zero, 
for otherwise the condition that $\lambda\cdot\delta=0$ 
implies that $\lambda_k=0$, and then
since $\lambda_j=0$ we have $\lambda=0$. 
Thus $\lambda_k = -\sum_{i\in I''} \delta_i \lambda_i \prec 0$.

By Lemma~\ref{l:dominantcase} there is some integer $r$ with
$(m\delta-\alpha')_i=r\delta_i$ for all $i$ with $\lambda_i\neq 0$.
Let $\beta = \alpha'-(m-r)\delta\in\Z^I$.
Of course $\beta_j=0$ and for any vertex $i$ with 
$\lambda_i\neq 0$ we have $\beta_i=0$.

Suppose that $\beta$ is nonzero.
Consider the restriction of $\beta$ to a connected
component of the quiver obtained from $Q'$ by
deleting all vertices $i$ with $\lambda_i\neq 0$.
It is actually a subquiver of $Q''$, so Dynkin.
If $\eta$ is a positive root for this connected
component, then $\eta\in R_\lambda^+$, and
\[
(\beta,\eta) = (\alpha',\eta) - (m-r)(\delta,\eta)
= (\alpha',\eta),
\]
so Lemma~\ref{l:rootineq} implies that
$-1 \le (\beta,\eta) \le 0$.
But this is impossible by Lemma~\ref{l:dynkvec} below. 

Thus $\beta=0$, so $\alpha=\epsilon_j+(m-r)\delta$.
Now since $p(\alpha)=0$ we have $m=r$.
\end{proof}

The above proof uses the following result about
Dynkin quivers.

\begin{lem}\label{l:dynkvec}
If $Q^\circ$ is a Dynkin quiver with vertex set $I^\circ$
then there is no nonzero vector $\alpha\in \Z^{I^\circ}$
with $-1 \le (\alpha,\eta) \le 0$ for all positive
roots $\eta$ for $Q^\circ$.
\end{lem}

\begin{proof}
We cannot have $(\alpha,\epsilon_i)=0$ for all $i$,
for otherwise $(\alpha,\alpha)=0$, so $\alpha=0$
since $Q^\circ$ is Dynkin.

Embed $Q^\circ$ in an extended Dynkin quiver of the same
type by adding an extending vertex $s$, and consider
$\alpha$ as a dimension vector for this quiver.
Let $\delta$ be the minimal positive imaginary root.

Since $\delta-\epsilon_s$ is a root for $Q^\circ$ we have
$(\alpha,\delta-\epsilon_s)\ge -1$. 
Now it is equal to $\sum_{i\neq s} \delta_i (\alpha,\epsilon_i)$,
and all terms in the sum are $\le 0$, but not all are zero.
Thus exactly one term
is nonzero, say for $i=r$, and it is equal to $-1$. 
This implies that $r$ is an extending vertex, 
and $(\alpha,\epsilon_r)=-1$.
Thus the vector $-\alpha$ and the extending vertices $r$ and $s$ 
contradict \cite[Lemma 8.8]{CBmm}.
\end{proof}

\section{Decomposing the quiver}\label{s:someisos}
In this section we suppose that $Q$ is a quiver whose
vertex set $I$ is a disjoint union $\J\cup\K$,
and we write any $\alpha\in\N^I$ as 
$\alpha=\alpha_{\J} + \alpha_{\K}$
where the summands have support in $\J$ and $\K$ respectively.

\begin{lem}\label{l:Nintotwo}
If the dimension vector of any composition factor
of a $\pr$-module of dimension $\alpha$ has support
contained either in $\J$ or in $\K$ then
\[
N(\lambda,\alpha) \cong N(\lambda,\alpha_{\J}) \times N(\lambda,\alpha_{\K}).
\]
\end{lem}

\begin{proof}
We can identify
\[
\mu_{\alpha_{\J}}^{-1}(\lambda)\times
\mu_{\alpha_{\K}}^{-1}(\lambda)
\]
with a $\GG(\alpha)$-stable closed subvariety of
$\mu_{\alpha}^{-1}(\lambda)$ (defined by the vanishing
of all arrows with one end in $\J$ and the other end in $\K$).
The inclusion thus induces a closed embedding
\[
N(\lambda,\alpha_{\J}) \times N(\lambda,\alpha_{\K})\to N(\lambda,\alpha),
\]
and by the assumption on composition factors this is a bijection.
\end{proof}

We give some cases when this can be applied.
First we need a lemma.

\begin{lem}\label{l:onesidemod}
Suppose there is a unique arrow with one end in $\J$ and 
the other in $\K$, say connecting vertices $j\in \J$ 
and $k\in \K$. 
Let $\tilde{Q}$ be the quiver with vertex set 
$\K\cup \{j\}$ containing this arrow, and
all arrows with head and tail in $\K$.
Let $\mu$ be the vector for $\tilde{Q}$  
whose restriction to $\K$ is the same as $\lambda$, 
and with $\mu_j=0$.

Let $\alpha\in\N^I$ and assume that $\alpha_j=1$
and $\lambda\cdot\alpha_{\J}=\lambda\cdot\alpha_\K=0$.
Then $\alpha\in \N R_\lambda^+$ if and only if
$\alpha_\J\in\N R_\lambda^+$ and 
$\epsilon_j+\alpha_\K\in\N R_\mu^+$.
\end{lem}

\begin{proof}
The statement does not depend on the orientation
of the arrows in $Q$, so we may suppose that the arrow 
connecting $\J$ and $\K$ is $b:k\to j$.

By \cite[Theorem 4.4]{CBmm} the 
condition that $\alpha\in \N R_\lambda^+$ is
that there is a $\pr$-module of dimension $\alpha$.
Similarly for the other two conditions.

Now if the module is given by
an element $x\in\Rep(\overline{Q},\alpha)$, then
for any vertex $i$ we have
\[
\sum_{h(a)=i} x_a x_{a^*}
-
\sum_{t(b)=i} x_{a^*} x_a
=
\lambda_i 1.
\]
Taking the trace and summing over all $i\in\J$, all but
one term cancels, leaving $\tr(x_b x_{b^*})= 0$.
Since this is a $1\times 1$ matrix we have $x_b x_{b^*}=0$.
It follows that the components of $x$ corresponding
to arrows with head and tail in $\J$ define
a $\pr$-module of dimension $\alpha_\J$,
and the remaining components of $x$ define
a $\Pi^\mu(\tilde{Q})$-module of dimension $\epsilon_j+\alpha_\K$.
Clearly two such modules can also be used to construct a
$\pr$-module of dimension $\alpha$.
\end{proof}

\begin{lem}\label{l:caseII} 
Suppose that $\lambda\cdot\alpha_{\J}=0$,
there is a unique arrow $b$ with one end in $\J$ and the other in
$\K$, say connecting vertices $j\in \J$ and $k\in \K$, and 
$\alpha_j = \alpha_k = 1$.
Then the dimension vector of any composition factor of a 
$\pr$-module of dimension $\alpha$ has support contained in $\J$ or $\K$.
\end{lem}

\begin{proof}
Because of the existence of a module of
dimension $\alpha$ we have $\lambda\cdot\alpha=0$,
hence also $\lambda\cdot\alpha_\K=0$.
For a contradiction, suppose there is a composition factor 
whose dimension $\beta$ does not have support in $\J$ or $\K$.
Then $\beta_j = \beta_k = 1$. Since the dimension vector $\gamma$
of any other composition factor must have support in $\J$ or $\K$, 
and has $\lambda\cdot\gamma=0$, we deduce that
$\lambda\cdot\beta_\J = \lambda\cdot\beta_\K=0$.

By Lemma \ref{l:onesidemod} we have 
$\beta_\J\in\N R_\lambda^+$, and by symmetry also
$\beta_\K\in\N R_\lambda^+$.
But clearly $(\beta_\K,\beta_\J)=-1$, so that
$p(\beta) = p(\beta_\J)+p(\beta_\K)$,
contradicting the fact that $\beta\in\Sigma_\lambda$.
\end{proof}

\begin{lem}\label{l:caseIII}
Suppose that $\lambda\cdot\alpha_{\J}=0$,
there is a unique arrow with one 
end in $\J$ and the other in $\K$, say connecting
vertices $j\in \J$ and $k\in \K$, $\alpha_j = 1$, the restriction
of $Q$ to $\K$ is extended Dynkin with extending vertex $k$ 
and minimal positive imaginary root $\delta$,
and $\alpha_\K = m\delta$ with $m\ge 2$.
Then the dimension vector of any composition factor of a 
$\pr$-module of dimension $\alpha$ has support contained in $\J$ or $\K$. 
\end{lem}

\begin{proof}
Because of the existence of a module of
dimension $\alpha$, we have $\lambda\cdot\alpha_\K=0$.
Since the field $K$ has characteristic zero, we deduce that 
$\lambda\cdot\delta=0$.

For a contradiction, suppose there is a composition factor 
whose dimension $\beta$ does not have support in $\J$ or $\K$.
Then $\beta_j = 1$. Since the dimension vector $\gamma$
of any other composition factor must have $\gamma_j=0$, 
it has support in $\J$ or $\K$, and since it
has $\lambda\cdot\gamma=0$, we deduce that
$\lambda\cdot\beta_\J = \lambda\cdot\beta_\K=0$.
Also $m\delta-\beta_\K\in \N R_\lambda^+$.

Let $\tilde{Q}$ be the quiver obtained from $Q$ 
as in Lemma~\ref{l:onesidemod}, and let $\mu$ be the
corresponding vector. 
Since $m\delta-\beta_\K$ has support in $\K$ it
can be considered as an element of $\N R_\mu^+$.
By Lemma~\ref{l:onesidemod} we have
$\beta_\J\in \N R_\lambda^+$ and $\epsilon_j+\beta_\K \in \N R_\mu^+$.
Now by assumption $\beta_\K$ is nonzero, 
so Proposition~\ref{l:maincase} implies
that $\epsilon_j+\beta_\K\notin\Sigma_\mu$.
By \cite[Theorem 5.6]{CBmm} this implies that there 
are nonzero $\phi,\psi\in\N R_\mu^+$ with
$\phi+\psi=\epsilon_j+\beta_\K$ and $(\phi,\psi)_{\tilde{Q}}\ge -1$.
Without loss of generality, $\phi_j=0$ and $\psi_j=1$.
Considered as a dimension vector for $Q$ we
clearly have $\phi\in\N R_\lambda^+$.
Also, Lemma~\ref{l:onesidemod} applies to the
dimension vector $\psi+\beta_\J-\epsilon_j$,
and shows that it belongs to $\N R_\lambda^+$.
Since also
\[
(\phi,\psi+\beta_\J-\epsilon_j) = (\phi,\psi)_{\tilde{Q}}\ge -1,
\]
we have $\beta = \phi + (\psi+\beta_\J-\epsilon_j) \notin \Sigma_\lambda$
by \cite[Theorem 5.6]{CBmm}. A contradiction.
\end{proof}

\section{Proof of the theorem}
\begin{proof}[Proof of Theorem \ref{t:mainthm}]
We prove this for all $Q$, $\lambda$ and $\alpha\in\N R_\lambda^+$ 
by induction on the maximum possible number of terms 
in an expression for $\alpha$ as a sum of elements of $R_\lambda^+$.
If $\alpha\in\Sigma_\lambda$ then the assertions are vacuous,
so assume that $\alpha\notin\Sigma_\lambda$.

By \cite[Lemma 5.2]{CBmm} and Lemma~\ref{l:reflectionhomeo}
we can always apply a sequence of admissible reflections to 
the pair $(\lambda,\alpha)$.
Let $F_\lambda$ be the set of \cite[Section 7]{CBmm}.
If $\alpha\notin F_\lambda$ then by applying a sequence of admissible 
reflections  to $(\lambda,\alpha)$ we may assume that there is a 
loopfree vertex $i$ with $\lambda_i=0$ and $(\alpha,\epsilon_i)>0$.
Clearly in any decomposition of $\alpha$ as a sum of 
elements of $\Sigma_\lambda$ one of the terms, say $\beta$, 
has $(\beta,\epsilon_i)>0$. But by \cite[Lemma 7.2]{CBmm} 
this implies that $\beta = \epsilon_i$. 
Now $\alpha-\epsilon_i\in\N R_\lambda^+$, and
by the inductive hypothesis the assertions hold for
$\alpha-\epsilon_i$. If the decomposition is
\[
\alpha-\epsilon_i = \sigma^{(1)} + \dots + \sigma^{(r)}
\]
then clearly
\[
\alpha = \epsilon_i + \sigma^{(1)} + \dots + \sigma^{(r)}
\]
is a suitable decomposition of $\alpha$.
Moreover, if we have
\[
N(\lambda,\alpha-\epsilon_i) \cong 
\prod_{t=1}^s S^{m_t} N(\lambda,\sigma^{(t)}),
\]
then since $N(\lambda,\epsilon_i)$ is just a point,
any term $S^m N(\lambda,\epsilon_i)$
if it occurs, can be removed, and replaced by
$S^{m+1} N(\lambda,\epsilon_i)$
without changing the product. Thus by Lemma~\ref{l:dropsimple}
we obtain the required expression for $N(\lambda,\alpha)$.

Thus we are reduced to the case when 
$\alpha\in F_\lambda\setminus\Sigma_\lambda$.
By applying a sequence of admissible reflections to the pair $(\lambda,\alpha)$,
and then passing to the support quiver of $\alpha$,
we may assume that one of the cases (I), (II) or (III) 
of \cite[Theorem 8.1]{CBmm} holds. We deal with each of
these in turn.

Case (I). Here $Q$ is extended Dynkin, $\lambda\cdot\delta=0$,
and $\alpha=m\delta$ for some $m\ge 2$. 
By Lemma \ref{l:deltasum} and Theorem~\ref{t:sympower}
the decomposition $\alpha=\delta+\cdots+\delta$
has the required properties.

Case (II). Here $Q$ decomposes as in
Lemma~\ref{l:caseII}. In the notation of Section~\ref{s:someisos}
we write $\alpha = \alpha_\J + \alpha_\K$.
Since $\alpha\in\N R_\lambda^+$ there is a $\pr$-module
of dimension $\alpha$. Since the dimension vector
of any composition factor has support in $\J$ or $\K$
we deduce that $\alpha_\J$ and $\alpha_\K$ are in $\N R_\lambda^+$.
By the inductive hypothesis the conclusions of the theorem hold
for $\alpha_\J$ and $\alpha_\K$. 
Adding together the decompositions of $\alpha_\J$
and $\alpha_\K$ we obtain a decomposition of $\alpha$.
Obviously, since $\alpha_\J$ and $\alpha_\K$ have
disjoint support, no summand occurs in both parts.
The result thus follows from 
Lemmas \ref{l:Nintotwo} and \ref{l:caseII}.

Case (III). Here $Q$ decomposes as in Lemma~\ref{l:caseIII}.
We write $\alpha = \alpha_\J + m\delta$. Again $\alpha_\J$
and $m\delta$ are in $\N R_\lambda^+$ and by the inductive
hypothesis the conclusions of the theorem hold for them.
This gives a decomposition of $\alpha$ which has
the required
properties by Lemmas \ref{l:Nintotwo} and \ref{l:caseIII}.
\end{proof}

\begin{proof}[Proof of Proposition \ref{p:rootobs}]
(1) If $\beta$ is a real root in 
$\Sigma_\lambda$, then $N(\lambda,\beta)$ is a point
by \cite[Corollary 1.4]{CBmm}.

(2) If $\beta$ is an isotropic imaginary root in
$\Sigma_\lambda$, then it is indivisible, for if
$\beta = r\gamma$ then $\gamma$ is a root,
it has $\lambda\cdot\gamma=0$ since the base field
$K$ has characteristic zero, and the decomposition
$\beta=\gamma+\dots+\gamma$ has 
$p(\beta)<p(\gamma)+\dots+p(\gamma)$, contrary
to the definition of $\Sigma_\lambda$.

By \cite[Theorem 5.8]{CBmm}, some sequence of 
admissible reflections sends the pair $(\lambda,\beta)$
to a pair $(\lambda',\beta')$ with $\beta'$ in the 
fundamental region. Since it is isotropic imaginary 
we have $(\beta',\epsilon_i)=0$ for any vertex $i$ 
in the support of $\beta'$. By \cite[\S 1.2]{Kac}
this implies that the support 
quiver $Q'$ of $\beta'$ is extended Dynkin and 
$\beta'=\delta$, its minimal positive imaginary root.

Finally $N(\lambda,\beta)\cong N(\lambda',\delta)$
by Lemma~\ref{l:reflectionhomeo}, and this is
a deformation of the Kleinian singularity of type $Q'$
by Kronheimer's work \cite{Kronheimer}. See for example
\cite[Section 8]{CBH}. 

(3) Suppose that $\beta$ is a non-isotropic imaginary root 
in $\Sigma_\lambda$ and $m\ge 2$.
If $F_\lambda$ is the set of \cite[Section 7]{CBmm},
then \cite[Lemma 7.4]{CBmm}
implies that $\beta\in F_\lambda$, and hence also 
$m\beta\in F_\lambda$. 
Now in \cite[Theorem 8.1]{CBmm}, case (I) cannot
occur since $m\beta$ is non-isotropic, and cases (II) and (III)
cannot occur since all components of $m\beta$ are divisible by
$m$. Thus $m\beta\in\Sigma_\lambda$.
\end{proof}

\end{document}